\documentclass[a4paper,11pt]{article}
\usepackage{amsmath}
\usepackage{amsfonts}
\usepackage{amssymb}
\usepackage{enumerate}
\usepackage{graphics}
\RequirePackage{geometry}
\newtheorem{theorem}{Theorem}
\newtheorem{lemma}[theorem]{Lemma}
\newtheorem{corollary}[theorem]{Corollary}
\newtheorem{proposition}[theorem]{Proposition}

\newcommand{\prob}{\mathbb{P}}
\newcommand{\qed}{\hfill $\square$}
\newcommand{\proof}{\textit{Proof: }}

\begin{document}
\title{The gaps between the sizes of large clusters in $2 D$ critical percolation}
\author{J. van den Berg\footnote{CWI and VU University Amsterdam} \, and
R. Conijn\footnote{VU University Amsterdam} \\
%\, (and D. Kiss ?)\footnote{CWI, Amsterdam; starting October 2013 at University of Cambridge, UK} \\
{\footnotesize email: J.van.den.Berg@cwi.nl; R.P.Conijn@vu.nl}
% and D.Kiss@cwi.nl}
}
\date{}

\maketitle

\begin{abstract}
Consider critical bond percolation on a large $2 n \times 2 n$ box on the square lattice. It is well-known 
(see \cite{BCKS99}) that the size (i.e.
number of vertices) of the largest open cluster is, with high probability, of order $n^2 \pi(n)$, where $\pi(n)$
denotes the probability that there is an open path from the center to the boundary of the box. The same result holds
for the second-largest cluster, the third largest cluster etcetera.

J\'arai \cite{J03} showed that the differences between the sizes of these clusters is, with high probability,
at least of order $\sqrt{n^2 \pi(n)}$. Although this bound was enough for his applications (to
incipient infinite clusters), he believed, but had no proof, that the differences are in fact of the same order as
the cluster sizes themselves, i.e. $n^2 \pi(n)$. Our main result is a proof that this is indeed the case.

%Our main result is a proof that the differences between the sizes are indeed (with high probability) of
%order $n^2 \pi(n)$. 
\end{abstract}

\medskip\noindent
\begin{center}
{\small {\it 2010 Mathematics Subject Classification.} 60K35. }
\end{center}

\section{Introduction and statement of main results}
For general background on percolation we refer to \cite{G99} and \cite{BR06}. 
We consider bond percolation on the square lattice with parameter $p$ equal to its critical value $p_c = 1/2$.
Let $\Lambda_{n} = [-n,n]^{2} \cap \mathbb{Z}^{2}$ be the $2 n \times 2 n$ box centered at $0$ and let
$\partial \Lambda_{n} = \Lambda_{n} \setminus \Lambda_{n-1}$ be the (inner) boundary of the box.
For each vertex $v \in \mathbb{Z}^{2}$, we write $\Lambda_{n}(v) = \Lambda_{n} + v$.
Further, the open cluster in $\Lambda_n$ of the vertex $v$ is denoted by $\mathcal{C}_{n}(v)$.
More precisely,
\[
 \mathcal{C}_{n}(v) := \{ u \in \Lambda_{n}: u \leftrightarrow v \textrm{ inside } \Lambda_{n} \},
\]
where '$u \leftrightarrow v$ inside $\Lambda_{n}$' means that there is an open path from $u$ to $v$ of which all  vertices
are in $\Lambda_{n}$.
We write $\pi(n)$ for the probability $\prob(O \leftrightarrow \partial \Lambda_{n})$, the probability that
there is an open path from $O$ to $\partial \Lambda_n$. Further, we write
\begin{equation} \label{sn-def}
 s(n) : = n^{2}\pi(n).
\end{equation}
By the {\em size} of a cluster we mean the number of vertices in the cluster.
Let, for $i =1, 2, \cdots$, $\mathcal{C}_{n}^{(i)}$ denote the $i$-th largest open cluster in $\Lambda_{n}$,
and let $|\mathcal{C}_{n}^{(i)}|$ denote its size.
(If two clusters have the same size, we order them in some deterministic way).

In \cite{BCKS99} it was proved that $|\mathcal{C}_{n}^{(1)}|$ is of
order $s(n)$.
In the later paper \cite{BCKS01} by the same authors it is shown that also
$|\mathcal{C}_{n}^{(2)}|$, $|\mathcal{C}_{n}^{(3)}|$ etcetera are of order $s(n)$.
They also proved an extension of this result for the case where the parameter $p$ is not
equal but close to $p_{c}$.

It was shown by
J\'{a}rai in \cite{J03} that for each $i$ the difference
$|\mathcal{C}_{n}^{(i)}| - |\mathcal{C}_{n}^{(i+1)}| \to \infty$ in probability as $n \to \infty$ .
In fact he showed that this
difference is at least of order $\sqrt{s(n)}$. He suggested that it should be of order $s(n)$, but did not have a proof.
In this paper we show that his conjecture is correct.
% and can be proved by modifications and refinements of his proof of the above mentioned
%weaker bound.
We became interested in such problems through our investigation of frozen-percolation processes.
Our main theorem is as follows.
\begin{theorem}\label{thm:GapsClusters}
 For all $k \in \mathbb{N}, \delta >0$, there exist $\varepsilon > 0, N \in \mathbb{N}$ such that
for all $n \ge N$:
\begin{equation}\label{eq:thm:GC}
 \prob\left(\exists i \leq k - 1 \, : \; |\mathcal{C}_{n}^{(i)}| - |\mathcal{C}_{n}^{(i+1)}| \le \varepsilon s(n)\right) < \delta.
\end{equation}
\end{theorem}

\medskip\noindent
\textbf{Remarks:} {\it (i) The analog of Theorem \ref{thm:GapsClusters} can be proved for site
and bond percolation on other common two dimensional lattices,
e.g. site percolation on the square or the triangular lattice. In this latter model (site percolation on the triangular lattice)
one of the last steps of the proof can be made a little bit shorter (see the Remark below the proof of Proposition
\ref{prop:GapsClustersDiam}). \\
(ii) The proof, which is given in Section \ref{section-proof-main}, follows the main line of
J\'arai's proof of the weaker bound: We divide the box $\Lambda_n$ in boxes of smaller length (denoted by $2 t$), and
condition on the configuration outside certain open circuits in these smaller boxes. Conditioned on
this information, the `contributions' (to the sizes of certain open clusters) from the interiors
of these circuits are independent random variables. This leads to a problem concerning the concentration
function of a sum of independent random variables, to which a general (`classical') theorem is applied.
The main difference with J\'arai's arguments is that we take $t$ proportional to $n$, with a
proportionality factor chosen as a suitable function of the `parameters' $k$ and $\delta$ in the theorem.
This makes the
arguments more powerful (and also somewhat more complicated). Moreover, the theorem on concentration functions we used
(see Theorem \ref{thm:ConcentrationFunction} below) is somewhat stronger than the one used in J\'arai's arguments.}

\vspace{0.3cm}Furthermore, with essentially the same argument we can show that the probability that there exists a cluster
with size in a given interval of length $\varepsilon s(n)$ goes to zero as $\varepsilon \to 0$ uniformly in $n$:
\begin{theorem}\label{thm:ClustInSmallInt}
 For all $x, \delta > 0$, there exists an $\varepsilon > 0$ such that, for all $n \in \mathbb{N}$:
\begin{equation}\label{eq:thm:CSI}
 \prob\left(\exists u \in \Lambda_{n}: xs(n) < |\mathcal{C}_{n}(u)| < (x + \varepsilon)s(n)\right) < \delta.
\end{equation}
\end{theorem}
%From this theorem it immediately follows that,
%for all $x, \delta > 0$, there exist $\varepsilon > 0$ such that for all $n \in \mathbb{N}$:
%\[
% \prob\left(xs(n) < |\mathcal{C}^{(1)}| < (x + \varepsilon)s(n)\right) < \delta.
%\]
This last theorem is in some sense complementary to the result in an earlier paper \cite{BC12}, where
we proved that, for any interval $(a,b)$, the probability that $|\mathcal{C}_{n}^{(1)}|/s(n) \in (a,b)$ is
bounded away from zero as $n \rightarrow \infty$.
% This theorem implies the same statement for the maximal cluster which is complementary to the result in an earlier paper \cite{BC12},
% where we proved that, for any interval $(a,b)$, the probability that $|\mathcal{C}_{n}^{(1)}|/s(n) \in (a,b)$ is bounded away from zero,
% uniformly in $n$.

\section{Notation and Preliminaries} \label{sect-prelim}

\subsection{Preliminaries}
First we need some more notation. For a cluster $\mathcal{C}_{n}(u)$ we define its (left-right) diameter by
\[
 \textrm{diam}(\mathcal{C}_{n}(u)) = \max_{v,w \in \mathcal{C}_{n}(u)} |v_{1} - w_{1}|.
\]
For a box $\Lambda_{n}$ we define the spanning cluster by
\begin{equation}\label{eq:def:SC}
 SC_{n} = \{ u \in \Lambda_{n}: u \leftrightarrow L(\Lambda_{n}) \textrm{ and } u \leftrightarrow R(\Lambda_{n})\},
\end{equation}
where $L(\Lambda_{n}) = \{-n\} \times [-n,n] \cap \mathbb{Z}^{2}$ and
$R(\Lambda_{n}) = \{n\} \times [-n,n] \cap \mathbb{Z}^{2}$.
We use the notation $A_{m,n}$ for the annulus $\Lambda_{n} \setminus \Lambda_{m}$
and, for a vertex $v \in \mathbb{Z}^{2}$, the notation $A_{m,n}(v)$ for $A_{m,n} + v$.

In our proof of Theorem \ref{thm:GapsClusters} and \ref{thm:ClustInSmallInt} we will use the
following results from the literature,
Theorems \ref{thm:ratioPi} - \ref{thm:ConcentrationFunction} below.
The first one is well known, see for example \cite{BCKS01}, \cite{BK85}.
\begin{theorem}\textbf{(\cite{BCKS01},\cite{BK85})}\label{thm:ratioPi}
 There exist constants $c_{1}, c_{2}, c_{3} > 0$, such that for all $m \le n$:
   \begin{equation*}
    c_{1}(\frac{n}{m})^{c_{2}} \le \frac{\pi(m)}{\pi(n)} \le  c_{3}(\frac{n}{m})^{\frac{1}{2}}.
   \end{equation*}
\end{theorem}
As we already mentioned in the introduction, the largest clusters in $\Lambda_{n}$
are of order $s(n)$. This is stated in the following result.
\begin{theorem}\textbf{(\cite{BCKS01} Thm. 3.1(i), 3.3, 3.6) }\label{thm:SizeIthLargestCluster}
 For all $i \in \mathbb{N}$,
\begin{equation}\label{eq:thm:SICi}
 \mathbb{E}[|\mathcal{C}_{n}^{(i)}|] \asymp s(n),
\end{equation}
and,
%for all $i \in \mathbb{N}$
 \begin{equation}\label{eq:thm:SICii}
  \liminf_{n \to \infty} \prob\left( \varepsilon < \frac{|\mathcal{C}_{n}^{(i)}|}{\mathbb{E}[|\mathcal{C}_{n}^{(i)}|]} < \frac{1}{\varepsilon}\right) \to 1 \qquad \textrm{as } \varepsilon \to 0.
 \end{equation}
\end{theorem}
In an earlier paper Borgs, Chayes, Kesten and Spencer showed exponential decay for the probability that
there exists a cluster with large volume, but a small diameter:
\begin{theorem}\textbf{(\cite{BCKS99} Remark (xiii))}\label{thm:LargeVolSmallDiam}
 There exist $C_{1}, C_{2} > 0$ such that for all $x > 0$, $\alpha \in (0,1]$ and $n \ge 4/\alpha$ we have
\begin{equation}\label{eq:thm:LVSD}
 \prob\left(\exists u \in \Lambda_{n}: |\mathcal{C}_{n}(u)| \ge xs(n); \; \mathrm{diam}(\mathcal{C}_{n}(u)) \le \alpha n\right)
 \le C_{1}\alpha^{-2}\exp{(-C_{2}x/\alpha)}.
\end{equation}
\end{theorem}
An easy consequence of Theorems \ref{thm:SizeIthLargestCluster} and \ref{thm:LargeVolSmallDiam} is the following.
\begin{corollary}\label{cor:LargestClustersLargeDiam}
 Let $k \in \mathbb{N}$. For all $\delta > 0$ there exist $\alpha > 0$ and $N \in \mathbb{N}$ such that for all $n \ge N$:
\begin{equation}\label{eq:cor:LCLD}
 \prob\left(\exists i \le k: \mathrm{diam}(\mathcal{C}_{n}^{(i)}) < \alpha n\right) < \delta.
\end{equation}
\end{corollary}

In \cite{J03} a version of Theorem \ref{thm:SizeIthLargestCluster} for the spanning cluster is given:
\begin{theorem}\textbf{(\cite{J03} Thm. 8)}\label{thm:SpanningCluster}
% There exist $C_{1}, C_{2} > 0$ such that $\forall n \in \mathbb{N}$:
 \begin{equation}\label{eq:thm:SCi}
  \mathbb{E}[|SC_{n}|] \asymp s(n);
 \end{equation}
 moreover,
 \begin{equation}\label{eq:thm:SCii}
  \lim_{\varepsilon \to 0} \inf_{n \in \mathbb{N}} \prob\left(\varepsilon < \frac{|SC_{n}|}{\mathbb{E}[|SC_{n}|]} < \frac{1}{\varepsilon} \; | \; SC_{n} \neq \emptyset\right) = 1.
 \end{equation}
\end{theorem}

In the proof of our main theorem we use the following inequality concerning the concentration function $Q(X,\lambda)$
of a random variable $X$, which is defined by
\begin{equation} \label{Q-def}
 Q(X,\lambda) = \sup_{x \in \mathbb{R}} \prob(x \le X \le x + \lambda),
\end{equation}
for $\lambda > 0$.
\begin{theorem}\textbf{(\cite{C65}; \cite{E68} (B))}\label{thm:ConcentrationFunction}
 Let $(X_{k})_{k \in \mathbb{N}}$ be a sequence of independent random variables, and $0 < \tilde{\lambda} \le \lambda$.
Let $a > 0$ and let $(b_k)_{ k \in \mathbb{N}}$ be a sequence of real numbers such that,
for all $k \in \mathbb{N}$,
\[
 \prob(X_{k} \le b_{k} - \frac{\tilde{\lambda}}{2}) \ge a,\qquad \prob(X_{k} \ge b_{k} + \frac{\tilde{\lambda}}{2}) \ge a.
\]
There exists a universal constant $C > 0$ such that, for all $m \in \mathbb{N}$
\[
 Q(S_{m},\lambda) \le \frac{C\lambda}{\tilde{\lambda} \sqrt{m\, a}},
\]
 where $S_{m} = X_{1} + X_{2} + \cdots + X_{m}$.
\end{theorem}

%\vspace{0.5cm}\noindent\textbf{Remark} Notice that J\'{a}rai used the somewhat weaker Kolmogorov-Rogozin inequality.

\subsection{Large clusters contain many good boxes}
In the proof of our main theorem we need the following lemma, which is essentially already in \cite{J03}.
First some definitions. Recall the notation $A_{m,n}$ in the beginning of Section \ref{sect-prelim}..
Let $t \in 3\mathbb{N}$. (Later we will choose a suitable value for $t$).
For any $i,j \in \mathbb{Z}$ we say that the box $\Lambda_{t}(2ti,2tj)$ is `good' if there is
an open circuit in the annulus $A_{\frac{2}{3}t,t}(2ti,2tj)$; in that case we denote
the widest open circuit in that annulus by $\gamma_{i,j}$.
(Although $\gamma_{i,j}$ depends on $t$, we omit that parameter
from the notation).
%For a cluster $\mathcal{C}_{n}(u)$, a box $\Lambda_{t}(2ti,2tj)$ is
%called 'good' if $\gamma_{i,j}$ exists and is part of the cluster $\mathcal{C}_{n}(u)$.
For each vertex $u$ we denote by $G_{t}(\mathcal{C}_{n}(u))$ the set of good boxes in $\Lambda_n$
of which the
corresponding $\gamma_{i,j}$ is contained in the open cluster of $u$. More precisely,
\begin{equation}\label{eq:def:SetGoodBoxes}
 G_{t}(\mathcal{C}_{n}(u)) = \{ (i,j):
\Lambda_{t}(2ti,2tj) \subset \Lambda_{n} \mbox{ is good }; \gamma_{i,j} \subset \mathcal{C}_{n}(u) \}.
\end{equation}

\begin{lemma}\label{lem:LargeDiamManyGood}
Let $\alpha > 0$. For any $\delta, \beta > 0$ there exist $\eta > 0$ and $N \in \mathbb{N}$ such that,
for all $n \ge N$ and $t \in (0,\eta n) \cap 3\mathbb{N}$
\begin{equation}\label{eq:lem:LDMG}
 \prob\left( \exists u \in \Lambda_{n}: \mathrm{diam}(\mathcal{C}_{n}(u)) \ge
\alpha n; \; |G_{t}(\mathcal{C}_{n}(u))| < \beta \right) < \delta.
\end{equation}
\end{lemma}

J\'{a}rai proved a somewhat stronger statement (see ~(3.15) in \cite{J03}), but we only need this 
weaker statement and give a (short) proof.

\vspace{0.5cm}\noindent\proof 
The $C_i$'s in this proof denote universal constants larger than $0$. Their existence is
important but their precise value does not matter for the proof.
First note that we can cover the box $\Lambda_{n}$ by at most 
\begin{equation} \label{eq-C1}
\frac{C_1}{\alpha^2}
\end{equation}
rectangles of width $\frac{1}{4}\alpha n$ and length  $\frac{1}{2}\alpha n$,
such that every cluster with diameter at least $\alpha n$ crosses at least one of these rectangles in the easy direction.
We consider one such rectangle, namely $Q_{0} := [0,\frac{1}{4}\alpha n] \times [0, \frac{1}{2}\alpha n]$.
(The precise choice of $Q_0$ doesn't matter for our purpose).
By RSW and the BK inequality we have that the probability that there are more than $C_2$
disjoint horizontal open crossings of $Q_0$ is less than

\begin{equation} \label{eq-C2}
\frac{\delta}{2} \,\, \frac{\alpha^2}{C_1}.
\end{equation}
%(Notice that we may assume that $[0,\frac{1}{4}\alpha n] \times [0, \alpha n]$ is contained in $\Lambda_{n}$.)
Let $R_{l}$ denote the $l$-th lowest open crossing of $Q_{0}$.
We claim that there exist $\eta \in (0, 1/2)$ and $N \in \mathbb{N}$ such that, for any $n \ge N$,  deterministic
crossing $r_{0}$ of $Q_0$, and $t \in (0,\eta n)$,
\begin{equation}\label{eq:prf:LDMG:statRlr0}
 \prob\left( |G_{t}(\mathcal{C}_{n}(r_{0}))| < \beta \; | \; R_{l} = r_{0} \right) <
\frac{\delta\alpha^{2}}{2C_1 C_2},
\end{equation}
where $\mathcal{C}_{n}(r_{0})$ denotes the open cluster which contains the crossing $r_{0}$.
From this claim we get (see \eqref{eq-C1} and \eqref{eq-C2}) that the l.h.s. of \eqref{eq:lem:LDMG} is less than

$$\frac{C_1}{\alpha^2} \left(\frac{\delta}{2} \, \frac{\alpha^2}{C_1} + C_2 \frac{\delta\alpha^{2}}{2C_1 C_2}\right) = \delta,$$
and the lemma follows.

It remains to prove the claim concerning the inequality \eqref{eq:prf:LDMG:statRlr0}:
%Conditioned on the $l$-th lowest crossing of $Q_{0}$,
%the configuration above this crossing is unexplored, and hence independent.
The objects defined below involve a parameter $i$. We will always assume that $i$ is such that the
corresponding object is contained in $[0, \frac{1}{4} \alpha n] \times [0, \alpha n]$.
Consider all rectangles of the form $A(i) := [2ti-t, 2ti+t] \times[0, \alpha n]$.
%where
%$i \in \{ 1, \cdots, \lfloor (\frac{1}{4}\alpha n - t)/(2t)\rfloor \}$.
For every $i$ we let $j(i)$ be the smallest integer $j$
for which the box $\Lambda_{t}(2ti,2t j)$ is located above $r_{0}$.
%with $\Lambda_{t}(2ti,2tj(i)) \cap r_{0} = \emptyset$.
Let $E(i)$ be the event that $\Lambda_{t}(2ti,2tj(i))$ is good
and $\gamma_{i,j(i)}$ is connected with $r_{0}$ inside $A(i)$.
The events $E(i)$ are conditionally  independent of each other (where we condition on the event $R_l = r_0$),
and, by RSW, each has probability larger than $C_3$.
Hence, when $\eta$ is small enough (that is, $n/t$ and thus the
number of events $E(i)$ is large enough), the probability that at most $\beta$ of the $E(i)$'s occur
is smaller than the r.h.s. of \eqref{eq:prf:LDMG:statRlr0}. This proves the claim and completes the proof of
Lemma \ref{lem:LargeDiamManyGood}.\qed

\vspace{0.5cm}\noindent\textbf{Remark:} {\em In one of the steps of J\'{a}rai's proof (see the lines below our statement of 
Lemma \ref{lem:LargeDiamManyGood}), he shows that with large probability the $l$-th
lowest crossing in $Q_{0}$ is contained in
$[0,\frac{1}{4}\alpha n] \times [0, \frac{1}{2}\alpha n(1-a)]$, for some constant $a < 1$. He used this to guarantee
that the good boxes obtained are inside $Q_{0}$. However, as the above arguments show, this (and hence the introduction
of the extra constant $a$) is not needed in our argument}.

\section{Proof of Theorems \ref{thm:GapsClusters} and \ref{thm:ClustInSmallInt}} \label{section-proof-main}

\subsection{Gaps between sizes of clusters with large diameter}

The following lemma will be used later to show that the conditions for Theorem \ref{thm:ConcentrationFunction} are
satisfied in our situation.
First we define, for each circuit $\gamma$,
int$(\gamma)$ as the interior of $\gamma$ (that is, the bounded connected component of
$\mathbb{R}^{2} \setminus \gamma$, where $\gamma$ is seen as subset of the plane), and
\begin{equation} \label{def-X}
X_{\gamma} := |\{ u \in \mathrm{int}(\gamma) \cap \mathbb{Z}^{2}: u \leftrightarrow \gamma \}|.
\end{equation}

\begin{lemma}\label{lem:RandomnessInCirc}
 There exist universal constants $\chi, \xi > 0$ and, for all $t \in 3\mathbb{N}$ and
for any circuit $\gamma$ in $A_{\frac{2}{3}t,t}$, a value $c(t,\gamma) \ge 0$ such that
\begin{eqnarray}
 \prob(X_{\gamma} \le c(t, \gamma)) & \ge & \chi;\label{eq:lem:RICi}\\
 \prob(X_{\gamma} \ge c(t, \gamma) + \xi s(t)) & \ge & \chi. \label{eq:lem:RICii}
\end{eqnarray}
\end{lemma}
\proof Fix some $a \in (0,\frac{1}{2})$. Define the random variable
$Z = |\{u \in A_{\frac{1}{3}t,t} \cap \mathrm{int}(\gamma): u \leftrightarrow \gamma \}|$.
Let $c(t,\gamma)$ be defined by
\[
 c(t,\gamma) = \min \{z \in \mathbb{N} \cup \{ 0 \}: \prob(Z \le z) > a \}.
\]
By RSW, the probability that there is a closed dual circuit in $A_{\frac{1}{3}t,\frac{2}{3}t}$ is larger
than some universal constant $C_1 > 0$. Moreover, if there is such a circuit, then $X_{\gamma} = Z$.
Hence, 
$\prob(X_{\gamma} \le c(t,\gamma))$ is larger than or equal to the probability that there is
such a circuit and that $Z \le c(t,\gamma)$. By the above and FKG this is larger than $C_1 a$.

To prove \eqref{eq:lem:RICii} recall the notation \eqref{eq:def:SC} and define the random variable
$Y = |SC_{\frac{1}{3}t}|$. Theorem \ref{thm:SpanningCluster} implies that there exist
constants $C_{2}, \xi > 0$ such that, for all $t$, we have $\prob(Y \ge \xi s(t)) > C_{2}$.
Let $E$ be the event that there is an open crossing in $\Lambda_{\frac{1}{3}t}$ from top to bottom and
that this crossing is connected to $\gamma$. On $E$ we have that $X_{\gamma} \ge Z + Y$, since
the spanning cluster is connected to $\gamma$. 
By RSW, $\prob(E)$ is larger than some universal constant $C_3$. Hence

\begin{equation}
\prob\left(X_{\gamma} \ge c(t,\gamma) + \xi s(t)\right) \,\, \geq \,\, \prob(E; Z \geq c(t,\gamma); Y \geq \xi s(t)) 
\,\,  \geq  \,\, C_3 \, (1-a) \, C_2, 
\end{equation}
where the last inequality uses FKG.
This proves Lemma \ref{lem:RandomnessInCirc}. \qed

\medskip
Now we prove the following proposition, from which, as we show in the next subsection, Theorem \ref{thm:GapsClusters}
follows almost immediately.
The set of clusters with diameter larger than $\alpha n$ is denoted by $\mathbf{C}_{\alpha,n}$.
More precisely,
\begin{equation}\label{eq:def:SetLargeClusters}
 \mathbf{C}_{\alpha,n} = \{\mathcal{C}_{n}(u): u \in \Lambda_{n}; \mathrm{diam}(\mathcal{C}_{n}(u)) \ge \alpha n\}.
\end{equation}
\begin{proposition}\label{prop:GapsClustersDiam} For all $\alpha, \delta > 0$ there exist $\varepsilon = \varepsilon(\alpha, \delta) >0, N = N(\alpha, \delta) \in \mathbb{N}$
such that, for all $n \ge N$
\begin{equation}\label{eq:prf:GC:StateOnlyDiam}
 \prob(\exists \mbox{ distinct } \mathcal{D}_{1}, \mathcal{D}_{2} \in \mathbf{C}_{\alpha,n}:
||\mathcal{D}_{1}| - |\mathcal{D}_{2}|| < \varepsilon s(n)) < \delta.
\end{equation}
\end{proposition}
\proof Let $\alpha, \delta > 0$ be given. By a standard RSW argument, the probability that
$|\mathbf{C}_{\alpha,n}| \geq 1$ is smaller than some constant $< 1$ which depends only on $\alpha$.
Hence, by the BK inequality we can choose a $\kappa = \kappa(\alpha, \delta) \in \mathbb{N}$
such that, for all n:
\begin{equation}\label{eq:prf:GC:FewLargeClusters}
 \prob(|\mathbf{C}_{\alpha,n}| > \kappa) < \frac{\delta}{3}.
\end{equation}
Let $\xi$ and $\chi$ as in Lemma \ref{lem:RandomnessInCirc} and $C$ as in Theorem \ref{thm:ConcentrationFunction}.
Take $\beta$ so large that \begin{equation}\label{beta-ineq}
 \frac{\xi}{2} \le
\frac{\delta \xi \sqrt{\chi}}{6C {\kappa\choose 2}} \cdot \sqrt{\beta}.
\end{equation}
(For the time being, this property of $\beta$ will play no role; it will become essential
at \eqref{eps-ineq} for a suitable choice of $\varepsilon$).
Let $\eta$ be as in Lemma \ref{lem:LargeDiamManyGood} (but with $\delta/3$ instead of $\delta$ in \eqref{eq:lem:LDMG}).
It is clear from that lemma that without loss of generality we may assume that 
\begin{equation}
\label{et-assump}
\eta < \frac{\alpha}{2}.
\end{equation}

\noindent
For each $n$ we take $t = t(n) = 3\lfloor \frac{1}{3}\eta n\rfloor$.
%Lemma \ref{lem:LargeDiamManyGood} implies that, for all sufficiently large $n$,
Hence, by the above choice of $\eta$ we have, for all sufficiently large $n$,
\begin{equation}\label{eq:prf:GC:ManyGood}
 \prob\left( \exists \mathcal{D} \in \mathbf{C}_{\alpha,n}: |G_{t}(\mathcal{D})| < \beta \right) < \frac{\delta}{3}.
\end{equation}
Denote by $W$ the event that there are at most $\kappa$ clusters in $\Lambda_{n}$ with diameter
at least $\alpha n$ and all these clusters have at least $\beta$ good boxes.
Note that the complement of $W$ is the union of the event in the l.h.s. of \eqref{eq:prf:GC:FewLargeClusters} and
the event in the l.h.s. of \eqref{eq:prf:GC:ManyGood}, and hence has probability smaller than $2 \delta / 3$.
Therefore, to prove
Proposition \ref{prop:GapsClustersDiam} it is sufficient to show that there exists
$\varepsilon >0$ such that for all sufficiently large $n$,
\begin{equation}\label{eq:prf:GC:StateOnW}
 \prob(W \cap \{ \exists \mbox{ distinct } \mathcal{D}_{1}, \mathcal{D}_{2} \in \mathbf{C}_{\alpha,n}:
||\mathcal{D}_{1}| - |\mathcal{D}_{2}|| < \varepsilon s(n) \}) < \frac{\delta}{3}.
\end{equation}

% all edges inside $\gamma$
% and the corresponding vertices, except the vertices which are part of $\gamma$.
%The configuration of the edges in int$(\gamma)$ is independent of the rest. 
We define (compare with \eqref{eq:def:SetGoodBoxes}) 
\[
 G_{t,n} = \{ (i,j) \in \mathbb{Z}^{2}: \Lambda_{t}(2ti,2tj) \subset \Lambda_{n} \mbox{ is good } \}.
\]
Recall that we denote the outermost open circuit in $A_{\frac{2}{3}t,t}(2ti,2tj)$ (if it exists) by $\gamma_{i,j}$.
Denote the configuration on the edges in the set

\begin{equation} \label{H-def}
 H: = [-n,n]^{2} \setminus \left( \bigcup_{(i,j) \in G_{t,n}} \mathrm{int}(\gamma_{i,j}) \right)
\end{equation}

\noindent
by $\omega_{H}$.

To estimate the l.h.s. of \eqref{eq:prf:GC:StateOnW} we condition
first on the $\gamma_{i,j}$'s
and the configuration $\omega_{H}$.
Therefore, let $\tilde{G}$ be an arbitrary set of vertices  $(i,j)$ with
$\Lambda_{t}(2ti,2tj) \subset \Lambda_{n}$, and let,
for each $(i,j) \in \tilde{G}$, $\tilde{\gamma}_{i,j}$ be a (deterministic) circuit
in $A_{\frac{2}{3}t,t}(2ti,2tj)$. Let $\tilde H$ be the analog of \eqref{H-def}, with $\gamma$ replaced by
$\tilde \gamma$ and let $\tilde \omega$ be a configuration on $\tilde H$.
We will consider the conditional distribution
$\prob(\cdot | G_{t,n} = \tilde{G};
\,\, \gamma_{i,j} = \tilde{\gamma}_{i,j}\, \forall (i,j) \in G_{t,n};\,\, \omega_{H} = \tilde{\omega})$.
% by $\bar{\prob}(\cdot)$.
Note that the information we condition on allows us to distinguish all the
clusters in $\mathbf{C}_{\alpha,n}$
and their good boxes. (Here we used that \eqref{et-assump} implies that no
cluster of $\mathbf{C}_{\alpha,n}$ fits entirely
in the interior of one of the above mentioned $\gamma_{i,j}$'s).
We may assume that $\tilde \omega$ is such that $W$ holds.
%Hence $\bar{\prob}(W) \in \{0,1 \}$.
Let $\mathcal{D}_{1}, \mathcal{D}_{2}$ be two open clusters in $\mathbf{C}_{\alpha,n}$ for
the configuration $\tilde \omega$. Their sizes
can be decomposed as follows:
\begin{eqnarray}\label{eq:prf:GC:SizeD1}
% |\mathcal{D}_{1}| = |\mathcal{D}_{1} \cap H| + \sum_{(i,j) \in G_{t}(\mathcal{D}_{1})} X_{\gamma_{i,j}},
 |\mathcal{D}_{1}| = a_1 + \sum_{(i,j) \in G_{t}(\mathcal{D}_{1})} X_{\tilde\gamma_{i,j}}, \\ \nonumber
 |\mathcal{D}_{2}| = a_2 + \sum_{(i,j) \in G_{t}(\mathcal{D}_{2})} X_{\tilde\gamma_{i,j}}, \\ \nonumber
\end{eqnarray}
where  $a_1 = |\mathcal{D}_{1} \cap H|$ and $a_2 = |\mathcal{D}_{2} \cap H|$, and the $X$ variables are as 
defined in \eqref{def-X}.
The terms $a_1$ and $a_2$ can be considered as `fixed' (namely, determined by $\tilde \omega$), and
the $X_{\tilde\gamma_{i,j}}$'s as independent random variables.
Therefore, and because there are at most ${\kappa \choose 2}$ choices for $\mathcal{D}_{1}$ and $\mathcal{D}_{2}$,
to prove \eqref{eq:prf:GC:StateOnW} it is enough to show that there exists $\varepsilon > 0$, which does
not depend on $a_{1}, a_{2}, G_{t}(\mathcal{D}_{1}), G_{t}(\mathcal{D}_{2})$ and the $\tilde\gamma_{i,j}$'s, such that

\begin{equation}\label{eq:prf:GC:changeWithJarai}
%  \lefteqn{\bar{\prob}(||\mathcal{D}_{1}| - |\mathcal{D}_{2}|| < \varepsilon s(n))}\nonumber\\
\prob\left(\Bigg\lvert\left(a_1 + \sum_{\left(i,j\right) \in G_{t}\left(\mathcal{D}_{1}\right)} X_{\tilde\gamma_{i,j}}\right) - 
\left( a_2 + \sum_{\left(i,j\right) \in G_{t}\left(\mathcal{D}_{2}\right)} X_{\tilde\gamma_{i,j}}\right) \Bigg\rvert
<\varepsilon s\left(n\right)\right)
< \frac{\delta}{3 {\kappa\choose 2}},
\end{equation}

\noindent
On the event $W$ we have that $|G_{t}(\mathcal{D}_{1})| \ge \beta$. So we can mark $\beta$ of the good boxes in
$G_{t}(\mathcal{D}_{1})$, and condition (in addition to the earlier mentioned information) also on
the values of $X_{\gamma_{i,j}}$ for the remaining good boxes in $G_{t}(\mathcal{D}_{1})$ and
all the good boxes in $G_{t}(\mathcal{D}_{2})$.
Hence it is enough to show that there exists an $\varepsilon >0$ such that
\begin{equation}\label{eq:prf:GC:AlmostEsseen}
 \prob\left( |\sum_{m=1}^{\beta} X_{\gamma_{m}} - b| < \varepsilon s(n) \right) < \frac{\delta}{3 {\kappa\choose 2}},
\end{equation}
uniformly in $b \in \mathbb{N}$, $\gamma_{1}, \cdots, \gamma_{\beta}$ and $G_{t}(\mathcal{D}_{1})$,
where the $\gamma_{m}$'s are circuits in distinct annuli $A_{\frac{2}{3}t,t}(2ti,2tj)$.
%Recall the definition of the concentration function just above Theorem \ref{thm:ConcentrationFunction}.
%Since the left hand side of \eqref{eq:prf:GC:AlmostEsseen} is bounded from above by the concentration function,
We will do this by application of Theorem \ref{thm:ConcentrationFunction}, where Lemma \ref{lem:RandomnessInCirc}
(and our choice \eqref{beta-ineq} for $\beta$) enables a suitable application of that theorem:

From \eqref{beta-ineq} it follows immediately that for all $n$ there is
an $\varepsilon(n)$ such that

\begin{equation} \label{eps-ineq}
 \frac{\xi}{2}\cdot \frac{s(t)}{s(n)} \le
\varepsilon(n) \le
\frac{\delta \xi \sqrt{\chi}}{6C {\kappa\choose 2}} \cdot \sqrt{\beta}\cdot
\frac{s(t)}{s(n)}.
\end{equation}

\noindent
By the lower bound in Theorem \ref{thm:ratioPi}, the l.h.s. of \eqref{eps-ineq}
is bounded away from $0$, uniformly in $n$. Hence, $\inf_n \varepsilon(n) > 0$. Take
$\varepsilon$ equal to this infimum. We get (with $Q$ as in \eqref{Q-def},

\begin{eqnarray} \label{ineq-final}
 \prob\left( |\sum_{m=1}^{\beta} X_{\gamma_{m}} - b| < \varepsilon s(n) \right) & \le &
Q(\sum_{m=1}^{\beta} X_{\gamma_{m}}, 2\varepsilon s(n)) \nonumber \\ 
 \, & \le & Q(\sum_{m=1}^{\beta} X_{\gamma_{m}}, 2\varepsilon(n) s(n)) \nonumber \\
 \, & \le & \frac{2C}{\xi\sqrt{\beta \chi}}\cdot \frac{s(n)}{s(t)}\cdot \varepsilon(n),
\end{eqnarray}
where in the last inequality we used Lemma \ref{lem:RandomnessInCirc} and applied Theorem \ref{thm:ConcentrationFunction}
(with $\tilde \lambda = \xi s(t)$, $a = \chi$, $m = \beta$ and
$\lambda = 2 \varepsilon(n) s(n)$). Note that the condition $\tilde \lambda \le \lambda$ in that
theorem is satisfied because 
$\xi s(t) \leq 2 \varepsilon(n) s(n)$ by the first inequality in \eqref{eps-ineq}. \\
Now, by the second inequality of \eqref{eps-ineq} we have that the r.h.s. of \eqref{ineq-final} is
at most $\frac{\delta}{3 {\kappa\choose 2}}$. This shows \eqref{eq:prf:GC:AlmostEsseen}
and completes the proof of Proposition \ref{prop:GapsClustersDiam}.\qed

\medskip\noindent
{\bf Remark} {\em In the case of site percolation on the triangular lattice we can, in equation \eqref{eps-ineq} and the line above it, skip
the introduction of $\varepsilon(n)$, and choose $\varepsilon$ itself such that it is (for all sufficiently large $n$) between the l.h.s. and r.h.s. of
\eqref{eps-ineq}. For that percolation model such $\varepsilon$ exists because (see \cite{GPS13}, Proposition 4.9
and the last part of the proof of Theorem
5.1 in that paper) $\pi(t) / \pi(n)$, and hence $s(t)/s(n)$, has a limit as $n \rightarrow \infty$
(with $t/n$ fixed)}.

\subsection{Proof of Theorem \ref{thm:GapsClusters}}
Let $\delta$ and $k$ be fixed. By Corollary \ref{cor:LargestClustersLargeDiam} we can choose
$\alpha = \alpha(\delta, k)$ and $N_{1} = N_{1}(\delta, k)$ such that, for all $n \ge N_{1}$
\[
 \prob(\exists i \le k: \mathrm{diam}(\mathcal{C}_{n}^{(i)}) < \alpha n) < \frac{\delta}{2}.
\]
Further, by Proposition \ref{prop:GapsClustersDiam} there is an $\varepsilon >0$ such 
that the probability that there are two clusters with diameter larger than $\alpha n$ of which
the sizes differ less than $\varepsilon s(n)$
is smaller than $\delta /2$. Hence the l.h.s. of \eqref{eq:thm:GC} is less than $\delta/2 + \delta/2$.
%Theorem \ref{thm:GapsClusters} follows.
\qed

\subsection{Proof of Theorem \ref{thm:ClustInSmallInt}}
Let $x$ and $\delta$ be given.
By Theorem \ref{thm:LargeVolSmallDiam} we can find an $\alpha$ such that
\begin{equation*}
 \prob\left(\exists u \in \Lambda_{n}: |\mathcal{C}_{n}(u)| \ge xs(n); \; \mathrm{diam}(\mathcal{C}_{n}(u)) \le \alpha n\right) < \frac{\delta}{2}.
\end{equation*}
Let $\mathbf{C}_{\alpha,n}$ be defined as in \eqref{eq:def:SetLargeClusters}. It is enough to show that
there exist $\varepsilon = \varepsilon(\alpha,\delta) > 0, N = N(\alpha,\delta) \in \mathbb{N}$ such that,
for all sufficiently large $n$,
\begin{equation}
 \prob(\exists \mathcal{D} \in \mathbf{C}_{\alpha,n}: |\mathcal{D} - xs(n)| < \varepsilon s(n)) < \frac{\delta}{2}.
\end{equation}
This can be proved in practically the same way as Proposition \ref{prop:GapsClustersDiam}. 
(And, in fact, a bit easier, because now we deal with single clusters instead of pairs of clusters.
In particular the factor ${\kappa \choose 2}$ is replaced by $\kappa$ in the proof.) \qed

%\bibliographystyle{acm}
%\bibliography{mybiblio}

\bigskip\noindent
{\large\bf Acknowledgment} We thank Demeter Kiss for valuable discussions and for comments
on a draft of this paper.

\end{document}